%% file: main_arxiv.tex
\documentclass[11pt]{article}

\usepackage[a4paper,margin=2.5cm]{geometry}

\usepackage{amsmath,amssymb,amsthm}
\usepackage{physics}
\usepackage{bm}
\usepackage{graphicx}
\usepackage{booktabs}
\usepackage{multirow}
\usepackage[hidelinks]{hyperref}
\usepackage{enumitem}
\usepackage{mathtools}
\usepackage{cite}
\usepackage{float}

\title{
Neural Architectures as Functional Priors in Physics-Informed Control Problems
}

\author{
Sonia Rubio Herranz$^{1}$ \;
Fernando Carlos López Hernández$^{2}$ \;
Antonio López Montes$^{3}$
}

\date{}

\begin{document}

\maketitle


\begin{center}

{\small
$^{1}$Department of Statistics and Operations Research (EIO),\\
Universidad Complutense de Madrid, 28040 Madrid, Spain
}

\vspace{0.3cm}

{\small
$^{2}$Department of Software Engineering and Artificial Intelligence,\\
Universidad Complutense de Madrid, 28040 Madrid, Spain
}

\vspace{0.3cm}

{\small
$^{3}$Department of Analysis and Applied Mathematics,\\
Universidad Complutense de Madrid, 28040 Madrid, Spain
}

\vspace{0.5cm}

{\small
\texttt{sorubi01@ucm.es}
\quad
\texttt{fclh@ucm.es}
\quad
\texttt{bantonio@mat.ucm.es}
}

\end{center}

\vspace{0.5cm}

\begin{abstract}

\input{sections/abstract}

\end{abstract}


\section{Introduction}

\input{sections/section1_introduction}


\section{Control Problem Formulations}

\input{sections/section2_control}


\section{Experimental Analysis of the Linear RLC System}

\input{sections/section3_rlc}


\section{Non-linear Duffing Dynamics}

\input{sections/section4_duffing}


\section{Conclusions and future work}

\input{sections/section5_conclusions}
\nocite{*}
\bibliographystyle{unsrt}

\input{sections/main_arxiv.bbl}
\end{document}

%% file: sections/abstract.tex

\setlength{\parindent}{0pt}
\noindent

In this work we investigate the role of neural architectures as implicit
functional priors in control problems governed by ordinary differential
equations. Rather than focusing on highly complex problems, our objective is to
investigate architecture-dependent effects in controlled dynamical systems
within the simplest physically interpretable settings possible. In particular, we study a controlled linear RLC electrical circuit and a nonlinear
Duffing-type dynamical system. Both systems are analyzed first through
classical optimal-control formulations and later through PINN-based approaches. We compare different combinations of multilayer perceptrons (MLPs) and
Fourier-based KAN-like architectures, and analyze their influence on the
resulting controls. The numerical experiments suggest that different
architectural choices systematically generate qualitatively distinct controls,
even under identical governing equations, loss functionals, initial and target
states, training parameters and physical constraints. Significant differences appear in the spectral structure, smoothness, energy distribution, and phase-space behavior of the learned solutions. A central observation of this work is the emergence of a functional specialization phenomenon when the neural architectures are allowed sufficient freedom to shape the structure of the learned controls. More specifically, in the systems considered here,
Fourier-based architectures tend to produce trajectories with richer
oscillatory content, whereas smoother low-frequency-biased architectures tend
to generate more regular and energetically efficient controls. This suggests
that different functional components of the control problem may be handled more
efficiently by different neural architectures, leading to an implicit
specialization between state representation and control generation.

%% file: sections/section1_introduction.tex
Since their introduction by Raissi, Perdikaris and Karniadakis
in \cite{raissi2019}, Physics-Informed Neural Networks (PINNs) have become an
increasingly important framework in the numerical treatment of differential
equations, inverse problems, and scientific control problems \cite{wang2021understanding,wang2020pinnfailures}. The original
PINN formulation was primarily based on standard multilayer perceptron (MLP)
architectures, which remain one of the most widely used neural
parameterizations in Scientific Machine Learning.

More recently, alternative neural architectures such as Kolmogorov--Arnold Networks (KANs)
and Fourier-based architectures have begun to attract attention in the context of physics-
informed neural networks \cite{liu2024kan}.

New benchmark studies comparing standard PINNs and KAN-based approaches have reported that KAN-inspired models may achieve more accurate solutions and improved gradient estimates in several physics-informed learning problems \cite{benchmark2026kanpinn}. Most existing studies, however, focus primarily on approximation accuracy, convergence behavior, or computational efficiency. At the same time, recent developments in deep learning theory have revealed the existence of architecture-dependent phenomena such as spectral bias, implicit regularization \cite{neyshabur2017implicit}, and frequency-dependent learning effects \cite{rahaman2019spectral,xu2020frequency,jacot2018ntk}. These observations suggest that neural architectures may influence not only optimization behavior and approximation accuracy, but also the functional structure of the learned solutions.

This observation motivates the central question of the present work: whether different neural architectures can induce systematically different control strategies under identical physical constraints. In the present work we investigate the role of neural architectures as implicit
functional priors in control problems governed by ordinary differential
equations. Rather than focusing on highly complex systems, our objective is to
analyze architecture-dependent effects within the simplest physically
interpretable settings possible. In particular, we study a controlled linear
RLC electrical circuit together with a nonlinear Duffing-type dynamical system.

Both systems are analyzed first through classical optimal-control formulations
and later through PINN-based approaches. In particular, we introduce a
separated architecture strategy in which the physical state and the control
variable are represented through different neural architectures. This
allows us to investigate independently the role played by the architecture in
the representation of the system dynamics and in the generation of the control
signal.

The numerical experiments performed in this work reveal the emergence of a
functional specialization phenomenon. More specifically, in the systems
considered here, Fourier-based architectures tend to produce trajectories with
richer oscillatory content, whereas smoother low-frequency-biased
architectures tend to generate more regular and energetically efficient
controls. These observations suggest that different functional components of
the control problem may be handled more efficiently by different neural
architectures, leading to an implicit specialization between state
representation and control generation.

These results are consistent with interpreting neural architectures as implicit variational biases influencing the family of physically admissible solutions selected during training.

The remainder of the paper is organized as follows. Section 2 introduces both the classical optimal-control formulation and the PINN-based framework employed throughout the study. Section 3 presents the numerical experiments for the linear RLC control problem and analyzes the resulting architecture-dependent effects. Section 4 extends the study to the nonlinear Duffing dynamics and investigates the persistence of these effects in a nonlinear regime. Finally, Section 5 discusses the interpretation of the results, their implications for Scientific Machine Learning control problems, and possible directions for future work.

%% file: sections/section2_control.tex
To investigate the influence of neural architectures on learned control strategies, we consider two controlled dynamical systems of increasing complexity. The first is a linear RLC electrical circuit, which provides a fully controllable benchmark with a well-understood spectral structure. The second is a nonlinear Duffing-type oscillator, which introduces state-dependent nonlinearities while preserving physical interpretability. For both systems, we first describe the corresponding classical optimal-control formulation and subsequently introduce the PINN-based framework used throughout the experiments.
\subsection{Classical Optimal Control Formulation}

We first consider a controlled series RLC circuit. Let \(q(t)\) denote the electric charge and let
\begin{equation}
i(t)=q'(t)
\end{equation}

be the current. The controlled dynamics are given by
\begin{equation}
Lq''(t)+Rq'(t)+\frac{1}{C}q(t)=u(t),
\label{eq:rlc_second_order}
\end{equation}
where \(L>0\) is the inductance, \(R\geq 0\) is the resistance, \(C>0\) is the capacitance, and \(u(t)\) represents the externally applied voltage.

Equation \eqref{eq:rlc_second_order} is one of the simplest physically interpretable controlled dynamical systems. It combines inertia-like behavior through the inductive term \(Lq''\), dissipation through the resistive term \(Rq'\), elastic restoring behavior through the capacitive term \(C^{-1}q\), and external actuation through the voltage \(u(t)\).

The natural frequency of the undamped system is
\begin{equation}
\omega_0=\frac{1}{\sqrt{LC}},
\end{equation}
while the damping coefficient is determined by the normalized quantity \(R/L\) --- which make the damping independent of the inductance scale. Thus, before introducing neural approximations, the system already contains an
intrinsic spectral structure associated with its natural frequency and damping.

Introducing the state variable
\[
x(t)=
\begin{pmatrix}
q(t)\\
i(t)
\end{pmatrix},
\]
equation \eqref{eq:rlc_second_order} can be written as the first-order controlled system
\begin{equation}
\dot{x}(t)=Ax(t)+Bu(t),
\label{eq:linear_state_space}
\end{equation}
where \(A\) characterizes the natural evolution of the circuit in
the absence of external forcing. The term \(-1/(LC)\) determines the
oscillatory behavior through the interaction between the inductor and the
capacitor, while \(-R/L\) represents the damping effect introduced by the
resistance. \(B\) specifies how the external control voltage
\(u(t)\) acts on the system:
\begin{equation}
A=
\begin{pmatrix}
0 & 1\\
-\dfrac{1}{LC} & -\dfrac{R}{L}
\end{pmatrix},
\qquad
B=
\begin{pmatrix}
0\\
\dfrac{1}{L}
\end{pmatrix}.
\end{equation}

We seek an admissible voltage control \(u(t)\) capable of steering the system
from a prescribed initial state \(x_0\) to a desired target state \(x_T\)
over the time interval \([0,T]\). Accordingly, we impose the boundary
conditions:
\begin{equation}
x(0)=x_0=
\begin{pmatrix}
q_0\\
i_0
\end{pmatrix},
\qquad
x(T)=x_T=
\begin{pmatrix}
q_T\\
i_T
\end{pmatrix}.
\label{eq:Initial_Final_Conditions} 
\end{equation}

Before analyzing the neural-control formulation, it is important to verify that
the underlying linear system is controllable in the classical control-theoretic
sense. According to Kalman's controllability criterion, the pair ((A,B)) is
controllable if the controllability matrix
\[
\mathcal{C} = (B \quad AB)
\]
has full rank, which guarantees that the system can be driven from any initial
state to any desired target state by a suitable control input \cite{zhou1996robust,brunton2022data}.

In the present case,
\begin{equation}
B=
\begin{pmatrix}
0\\
\dfrac{1}{L}
\end{pmatrix},
\qquad
AB=
\begin{pmatrix}
\dfrac{1}{L}\\
-\dfrac{R}{L^2}
\end{pmatrix},
\end{equation}
which yields the controllability matrix
\[
\mathcal C=(B \quad AB).
\]
Since
\begin{equation}
\det(\mathcal{C})
=
-\frac{1}{L^2}
\neq 0,
\end{equation}
the matrix \(\mathcal C\) has full rank and  for every \(L>0\), the system is controllable. This means that, at least in the linear setting, the lack of controllability is not an obstruction because there exist controls capable of steering the circuit between arbitrary states in finite time.
This property makes the RLC model particularly suited for studying
architecture-dependent effects in neural control problems.

Among all admissible controls steering the system from \(x_0\) to \(x_T\), a natural classical choice is the minimum-energy control:
\begin{equation}
\min_{u}
J(u)
=
\int_0^T u(t)^2\,dt,
\label{eq:min_energy}
\end{equation}
subject to the controlled dynamics \eqref{eq:linear_state_space}
together with the initial and terminal conditions
\eqref{eq:Initial_Final_Conditions}.

For linear controllable systems, the minimum-energy control admits a closed-form
solution in terms of the finite-time controllability Gramian, which quantifies
how effectively the control input can influence the system over the interval
([0,T]) \cite{zhou1996robust}:

\begin{equation}
W_T
=
\int_0^T
e^{A(T-s)}BB^\top e^{A^\top(T-s)}
\,ds.
\end{equation}
If \(W_T\) is non-singular, the minimum-energy control is given by
\begin{equation}
u^*(t)
=
B^\top e^{A^\top(T-t)}
W_T^{-1}
\left(
x_T-e^{AT}x_0
\right),
\label{eq:min_energy_control}
\end{equation}

and the corresponding trajectory is obtained from the variation-of-constants formula
\begin{equation}
x(t)
=
e^{At}x_0
+
\int_0^t
e^{A(t-s)}B u^*(s)\,ds.
\end{equation}

\subsection{PINN-Based Control Formulation}

In order to approximate the solution of the controlled system \eqref{eq:linear_state_space} using neural-network representations 
\(\mathcal{N}(t)\), we employ separate neural networks for the state variables

\begin{equation}
(q(t),i(t))
=
\mathcal{N}_{\mathrm{state}}(t),
\end{equation}

and the control signal

\begin{equation}
u(t)
=
\mathcal{N}_{\mathrm{control}}(t).
\end{equation}

This separation allows us to investigate independently the influence of the architecture associated with the physical state and the architecture associated with the control variable.

\subsubsection*{Residual Equations}

The residuals measure the extent to which the neural-network approximations
violate the governing equations. They are defined as

\begin{equation}
r_1(t)=q'(t)-i(t),
\end{equation}

and

\begin{equation}
r_2(t)
=
i'(t)
+\frac{1}{LC}q(t)
+\frac{R}{L}i(t)
-\frac{1}{L}u(t).
\end{equation}

\subsubsection*{Loss Function}

A physically meaningful solution should simultaneously satisfy the governing
differential equations, the prescribed boundary conditions, and basic
regularity requirements on the control signal. To enforce these objectives
during training, the total loss functional is defined as

\begin{equation}
\mathcal{L}
=
\mathcal{L}_{\mathrm{ODE}}
+
\mathcal{L}_{\mathrm{IC}}
+
\mathcal{L}_{\mathrm{TC}}
+
\lambda_u
\int_0^T u(t)^2dt
+
\lambda_s
\int_0^T |u'(t)|^2dt.
\label{eq:total_loss2}
\end{equation}

The term \(\mathcal L_{\mathrm{ODE}}\) enforces the governing differential
equations by penalizing the residuals at \(N\) collocation points
\(\{t_j\}_{j=1}^{N}\subset(0,T)\):

\begin{equation}
\mathcal{L}_{\mathrm{ODE}}
=
\frac{1}{N}
\sum_{j=1}^{N}
\left(
|r_1(t_j)|^2
+
|r_2(t_j)|^2
\right).
\end{equation}

The terms \(\mathcal{L}_{\mathrm{IC}}\) and
\(\mathcal{L}_{\mathrm{TC}}\) enforce the prescribed initial and target
states through

\begin{equation}
\mathcal{L}_{\mathrm{IC}}
=
|q(0)-q_0|^2
+
|i(0)-i_0|^2,
\end{equation}

and

\begin{equation}
\mathcal{L}_{\mathrm{TC}}
=
|q(T)-q_T|^2
+
|i(T)-i_T|^2.
\end{equation}

The last two terms regularize the control signal.
The \(L^2\)-penalty

\[
\lambda_u
\int_0^T u(t)^2dt
\]

penalizes high-amplitude controls and therefore reduces the overall control
energy. On the other hand, the smoothness regularization

\[
\lambda_s
\int_0^T |u'(t)|^2dt
\]

discourages rapidly oscillating controls and suppresses high-frequency
components in the learned signal. The parameters
\(\lambda_u\geq0\) and \(\lambda_s\geq0\) determine the relative
importance of control-energy minimization and control smoothness.

\subsubsection*{Training Parameters}

In our experiments, we employ the following parameters in the loss
functional \eqref{eq:total_loss2}
\begin{equation} \lambda_{\mathrm{ODE}}=1, \qquad \lambda_{\mathrm{IC}}=100, \qquad \lambda_{\mathrm{TC}}=100, \end{equation} together with \begin{equation} \lambda_u=10^{-3}, \qquad \lambda_s=10^{-4}. \end{equation}

The weights were selected empirically after preliminary experiments aimed at
balancing physical accuracy and control regularity. The values
\(\lambda_{\mathrm{ODE}}=1\),
\(\lambda_{\mathrm{IC}}=\lambda_{\mathrm{TC}}=100\)
prioritize the satisfaction of the governing equations and boundary
conditions, while the relatively small regularization weights
\(\lambda_u=10^{-3}\) and \(\lambda_s=10^{-4}\) prevent unrealistic control
amplitudes and excessive oscillations without dominating the optimization
process. This choice preserves sufficient flexibility for the neural
architectures to influence the functional structure of the learned controls,
which is one of the central aspects investigated in this work. In addition, we employ \(N=900\) collocation points sampled uniformly over \([0,T]\). Preliminary sensitivity analyses indicated that increasing the number of collocation points beyond \(N=900\) produced only marginal changes in the
computed solutions while significantly increasing the computational cost.

This formulation defines the general PINN-based control framework employed in the present study. 

\subsection{Experimental Setup and Evaluation Metrics}

The experiments were implemented in Python using the PyTorch framework\footnote{Source code publicly available at \url{https://github.com/ifernandolopez/architecture_dependent_functional_priors}}.
The optimization process was carried out over 15000 training epochs
using the Adam optimizer with the same optimization parameters and
learning rate (\(5\times10^{-4}\)) across all architectural
configurations. Automatic differentiation was employed for the
computation of the residual terms appearing in the PINN loss
functional.

\subsubsection*{Neural Architectures}

We consider different combinations of MLP and FourierKAN
architectures for the state and control networks.
The detailed characteristics of the neural architectures employed in the experiments are summarized in Table 1.

\begin{table}[H]
\centering

\label{tab:architectures}

\begin{tabular}{lcccccc}
\toprule
Architecture &
Input &
Hidden Layers &
Width &
Output &
Representation &
Parameters \\
\midrule

MLP State
&
1
&
3
&
64
&
2
&
MLP + tanh
&
8578
\\

MLP Control
&
1
&
3
&
64
&
1
&
MLP + tanh
&
8513
\\

FourierKAN State
&
33
&
2
&
64
&
2
&
Fourier features
&
6466
\\

FourierKAN Control
&
33
&
2
&
64
&
1
&
Fourier features
&
6401
\\

\bottomrule
\end{tabular}

\caption{Neural architectures used in the experiments.}

\end{table}

For the Fourier-based architectures, the input dimension increases due to the sinusoidal feature expansion. In particular, we employ 16 Fourier modes, producing an effective input dimension equal to
\(
2\times16+1=33.
\)

\subsubsection*{Evaluation Metrics}

Our purpose is not merely to evaluate approximation accuracy, but to analyze how different neural architectures influence the physical structure of the learned controls. The evaluation therefore focuses simultaneously on dynamical accuracy, energetic efficiency, regularity, and spectral properties. The first quantity considered is the final-state accuracy. Since the objective of the control problem is to steer the system toward a prescribed target state, we measure the terminal errors
\[
|q(T)-q_T|,
\qquad
|i(T)-i_T|.
\]
These quantities indicate how accurately the learned control reaches the desired final configuration.

A second important metric is the residual error associated with the governing equations. Given the residual functions
\[
r_1(t)=q'(t)-i(t),
\]
and
\[
r_2(t)
=
i'(t)
+\frac{1}{LC}q(t)
+\frac{R}{L}i(t)
-\frac{1}{L}u(t),
\]
we evaluate the mean residual error through
\[
\|r\|^2
=
\frac{1}{N}
\sum_{j=1}^{N}
\left(
|r_1(t_j)|^2
+
|r_2(t_j)|^2
\right).
\]

The energetic cost of the control is evaluated through the quadratic energy functional
\[
E(u)
=
\int_0^T u(t)^2dt.
\]

In addition to the control energy, we also evaluate the temporal smoothness of the control through the Sobolev-type seminorm
\[
S(u)
=
\int_0^T |u'(t)|^2dt.
\]

Finally, we analyze the spectral structure of the learned controls using Fourier analysis. Given a control signal \(u(t)\), we compute its discrete Fourier transform and study the distribution of its frequency content. 

This spectral analysis plays a central role in the present work.

Taken together, these metrics provide a multi-perspective characterization of the learned
controls.

%% file: sections/section3_rlc.tex
\subsection{Control Results and Phase-Space Dynamics}

In this experiment, we study an initially charged capacitor with zero initial current. The objective of the control is therefore to drive the circuit toward the equilibrium configuration while compensating for the intrinsic oscillatory and dissipative dynamics of the system.

\subsubsection*{Experimental Configuration}

We consider the controlled RLC system with parameters
\[
L = 1, \qquad R = 0.4, \qquad C = 1.
\]

For these values, the circuit operates in an underdamped regime, allowing
oscillatory dynamics to remain visible while introducing moderate dissipation
through the resistive term. This parameter choice is intended to provide a
representative and physically interpretable test case rather than an exhaustive
parametric study of the RLC model.

The undamped circuit has natural frequency
\[
\omega_0=\frac{1}{\sqrt{LC}}=1,
\]
corresponding to the natural period
\[
T_0 = 2\pi.
\]
The control horizon is chosen as
\[
T = 8,
\]
which is slightly larger than one natural oscillation period. This allows the intrinsic oscillatory dynamics of the circuit to develop while keeping the control problem within a finite and nontrivial time scale.

The system is steered from the initial state
\[
(q_0,i_0)=(1,0)
\]
toward the target equilibrium state
\[
(q_T,i_T)=(0,0).
\]

Similar qualitative behaviors were observed across several nearby parameter
configurations, including variations of the damping coefficient and the control
horizon. These additional tests are not reported in detail, since the objective
of the present section is to isolate architecture-dependent effects in a fixed
and interpretable dynamical regime rather than to provide a complete
sensitivity analysis with respect to the physical parameters.

\begin{figure}[H]
\centering

\includegraphics[width=0.85\textwidth]{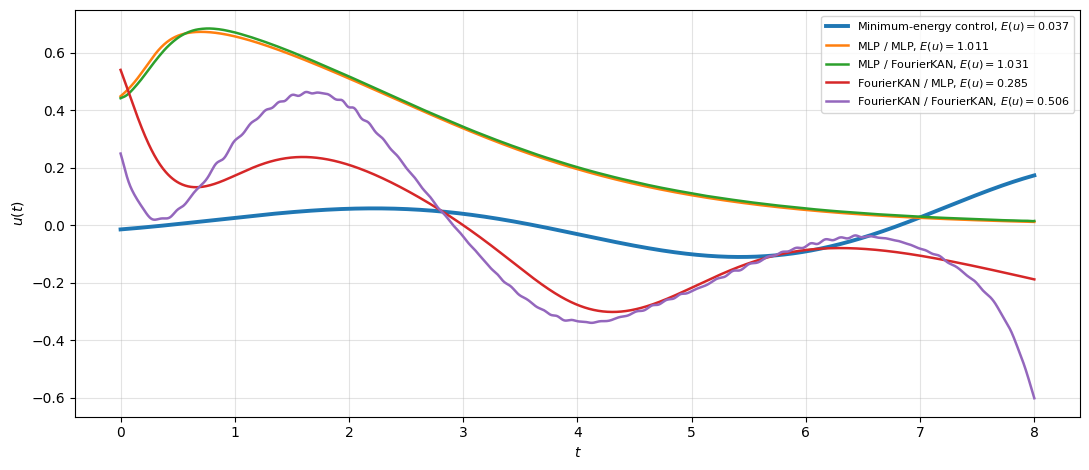}

\caption{Learned controls obtained with the different neural architectures together with the classical minimum-energy control.}

\label{fig:training_loss}

\end{figure}

\subsubsection*{Learned Control Signals}

Figure 1 shows the controls obtained with the different neural architectures together with the classical minimum-energy control. Although all configurations are trained under identical physical constraints and optimization settings, they generate qualitatively different control profiles. Differences in smoothness, oscillatory behavior, and overall functional structure are already visible, suggesting that the neural architecture influences the family of admissible controls selected during training.

\subsubsection*{State Trajectories and Phase-Space Dynamics}

Figure 2 shows the corresponding phase-space trajectories associated with the learned controls. Although all configurations successfully steer the system toward the same target state, noticeable differences appear in the geometry of the resulting trajectories. These differences indicate that the architectures are not merely producing alternative control signals, but are inducing distinct dynamical evolutions of the underlying physical system. In particular, the Fourier-based architectures tend to generate richer oscillatory motions in phase space, whereas the MLP-based configurations follow smoother and more regular paths. This observation provides further evidence that the neural architecture influences the family of admissible control strategies selected during training.

\begin{figure}[H]
\centering

\includegraphics[width=0.85\textwidth]{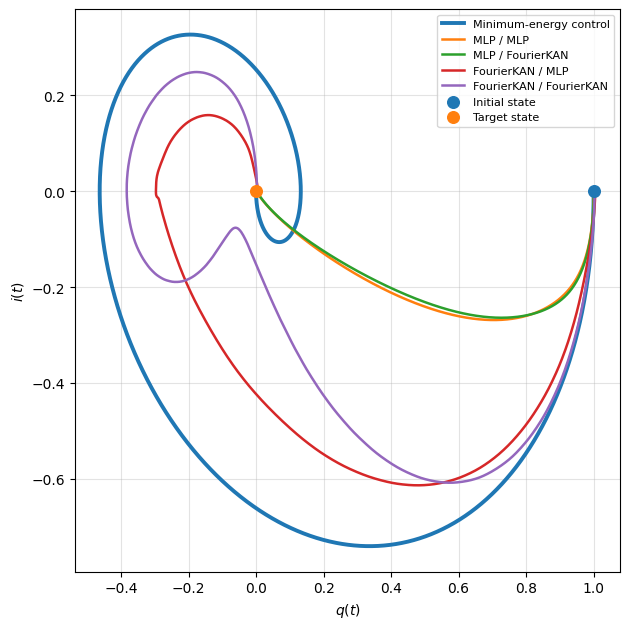}

\caption{
Phase-space trajectories corresponding to the different control strategies.
}

\end{figure}

\begin{figure}[H]
\centering

\includegraphics[width=0.85\textwidth]{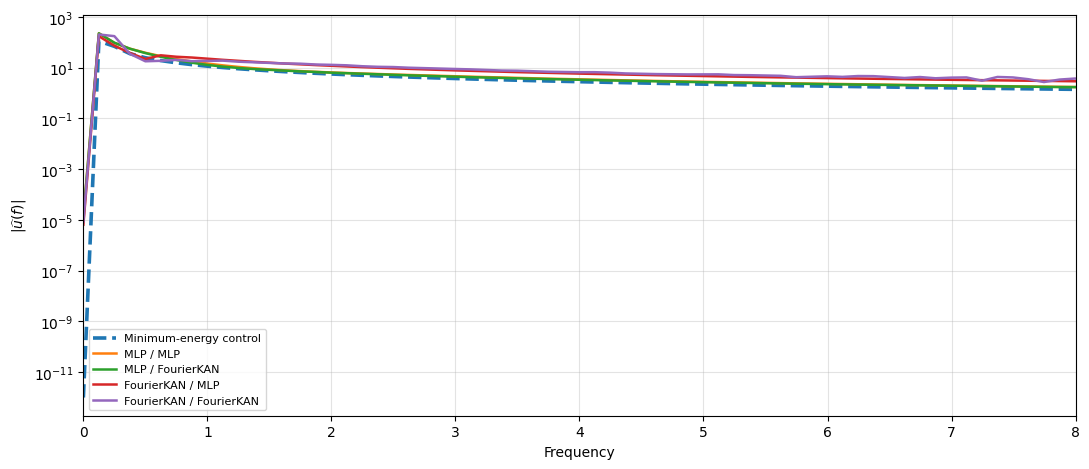}

\caption{
Fourier spectra associated with the different control strategies.
}

\end{figure}

Figure 3 shows the Fourier spectra associated with the learned controls in
logarithmic scale. The spectral distributions reveal clear
architecture-dependent differences in the frequency content of the learned
signals. In particular, the Fourier-based configurations preserve stronger
high-frequency components and exhibit a broader spectral distribution, whereas
the MLP-based controls concentrate most of their energy in the low-frequency
range. These observations are consistent with the spectral-bias phenomenon and
provide quantitative evidence that different neural architectures favor
different families of admissible controls. This tendency is further
confirmed quantitatively by the spectral-centroid values reported in Table 2.

\subsection{Spectral Analysis and Architecture-Dependent Effects}

Although all neural configurations successfully steer the system toward the prescribed target state, the resulting controls exhibit substantially different spectral and dynamical characteristics. These differences are visible both in the Fourier spectra of the learned controls and in the geometry of the corresponding phase-space trajectories.

The Fourier spectra reveal a systematic architecture-dependent behavior. While the MLP-based configurations concentrate most of their energy in the low-frequency range, the Fourier-based architectures preserve a richer high-frequency content. Although such tendencies are broadly consistent with the spectral-bias phenomenon reported in deep learning theory \cite{rahaman2019spectral}, our contribution is not merely to observe spectral differences between architectures. Rather, we show that these architectural preferences persist within a physics-informed optimal control framework and lead to distinct admissible control laws and dynamical trajectories.

\begin{table}[H]
\centering
\small
\begin{tabular}{|l|c|c|c|c|c|}
\hline
\textbf{Control Strategy} &
\textbf{$E(u)$} &
\textbf{Smoothness} &
\textbf{Spectral Centroid} &
\textbf{$q(T)$} &
\textbf{$i(T)$} \\
\hline

Minimum-energy control
& 0.037370
& --
& --
& -0.000004
& -0.000232 \\
\hline

MLP / MLP
& 1.011405
& 0.181780
& 12.654719
& -0.000649
& 0.001344 \\
\hline

MLP / FourierKAN
& 1.036404
& 0.186792
& 12.632654
& 0.001439
& 0.000236 \\
\hline

FourierKAN / MLP
& 0.285213
& 0.255568
& 16.579145
& 0.002404
& 0.000351 \\
\hline

FourierKAN / FourierKAN
& 0.508416
& 1.334480
& 16.760155
& 0.000628
& 0.000636 \\
\hline

\end{tabular}

\caption{
Summary of the main quantitative metrics associated with the different control strategies.
}

\label{tab:rlc_summary}
\end{table}

The quantitative metrics reported in Table 2 confirm these observations. The Fourier-based configurations exhibit larger spectral centroids, indicating a systematic shift toward higher frequencies. This systematic increase is accompanied by larger roughness values, reflecting the richer oscillatory structure already observed in Figures 1--3.

At the same time, all neural architectures achieve very small terminal errors in both \(q(T)\) and \(i(T)\), demonstrating that the observed differences are not caused by a failure to solve the control problem. Rather, they correspond to different admissible strategies for achieving the same dynamical objective.

%% file: sections/section4_duffing.tex
The results obtained for the linear RLC system suggest that neural
architectures may act as implicit functional priors, leading to
systematically different control strategies under identical physical
constraints. An important question is whether these architecture-dependent
effects are specific to linear dynamics or persist in the presence of
nonlinear phenomena.

To investigate this issue, we extend the analysis to a nonlinear
Duffing-type oscillator. Unlike the RLC circuit, the Duffing system
contains a cubic restoring term that introduces state-dependent
nonlinearities and destroys the superposition principle. This provides a
natural testbed for evaluating the robustness of the observed functional
specialization phenomenon beyond the linear regime.

\subsection{Duffing Control Problem}

We consider the nonlinear Duffing-type equation

\begin{equation}
q''(t)+\delta q'(t)+\alpha q(t)+\beta q(t)^3 = u(t),
\label{eq:duffing_controlled}
\end{equation}
where \(q(t)\) denotes the state variable, \(q'(t)\) its velocity-like component,
\(\delta\geq 0\) is a damping coefficient, \(\alpha\) is the linear stiffness
parameter, \(\beta\) measures the strength of the cubic nonlinearity, and
\(u(t)\) represents the external control.

The Duffing equation is one of the classical prototypical models in nonlinear
dynamics. It was originally introduced in the study of nonlinear vibrations and
has since appeared in a wide range of applications, including nonlinear
electrical circuits, mechanical oscillators, structural vibrations,
microelectromechanical systems, and nonlinear resonance phenomena
\cite{strogatz2018nonlinear}.

Unlike linear oscillators, the Duffing system contains the cubic restoring term \( \beta q(t)^3 \), which introduces an amplitude-dependent stiffness. As a consequence, the dynamical response depends on the current state amplitude and the superposition principle no longer holds. This nonlinear behavior gives rise to richer dynamical structures and more complex control landscapes than those encountered in linear systems.

For these reasons, the Duffing oscillator provides a natural benchmark for assessing whether the architecture-dependent functional specialization observed in the linear RLC problem remains present under nonlinear dynamics.

\subsection{Classical Baseline and PINN-Based Control Results}

Unlike the linear RLC system, the nonlinear Duffing equation does not admit an explicit minimum-energy control formula in terms of controllability Gramians. Consequently, the classical reference solution was obtained numerically through direct optimization of a discretized control signal.

More precisely, the temporal interval \( [0,T] \) was discretized into \( N \) time steps,
\begin{equation}
0=t_0<t_1<\cdots<t_{N}=T,
\end{equation}
and the control was represented by a discretized control vector
\begin{equation}
u=(u_0,u_1,\ldots,u_{N}).
\end{equation}

For each candidate control, the Duffing dynamics
\begin{equation}
q''(t)+\delta q'(t)+\alpha q(t)+\beta q(t)^3=u(t)
\end{equation}
were numerically integrated using a standard time-stepping scheme. This procedure produces the corresponding state trajectory and its final state.

The objective functional was defined as follows, where the first term measures the control energy, while the second penalizes deviations from the prescribed target state.

\begin{equation}
J(u)=\int_0^T u(t)^2dt
+\lambda_T |x(T)-x_T|^2,
\end{equation}

The optimization problem was then solved using gradient-based minimization over the discretized control variables. 
This numerical solution provides a classical nonlinear control baseline against which the different PINN-based architectures can be compared.

\subsection*{Learned Control Signals}

Figure 4 displays the controls obtained with the different neural architectures together with the classical nonlinear baseline. Although all configurations are trained under the same physical constraints and optimization objectives, they generate qualitatively different control profiles.

The differences observed in the linear RLC problem persist in the nonlinear setting. Variations in smoothness, oscillatory behavior, and overall functional structure are clearly visible across architectures. These differences suggest that the neural architecture continues to influence the family of admissible controls selected during training, even after the introduction of nonlinear dynamics.

\begin{figure}[H]
\centering

\includegraphics[width=0.85\textwidth]{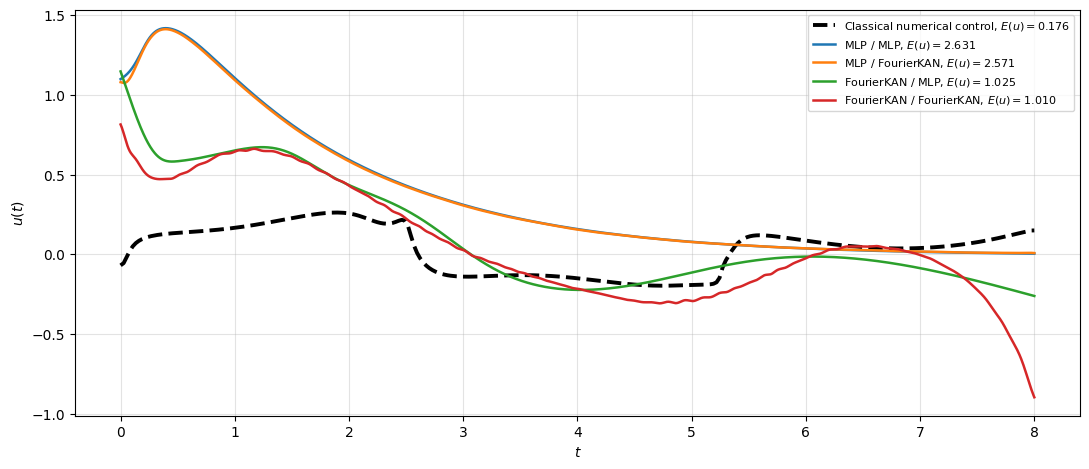}

\caption{
Control shapes.
}

\end{figure}

\subsection*{State Trajectories and Phase-Space Dynamics}

Figure 5 shows the phase-space trajectories associated with the different control strategies in the nonlinear Duffing setting. Although all configurations successfully steer the system toward the prescribed target state, substantial differences appear in the geometry of the resulting trajectories.

These differences indicate that the architectures are not merely generating alternative control signals. Rather, they induce distinct dynamical evolutions of the underlying nonlinear system. The trajectories follow different paths through phase space before reaching the same final configuration, revealing that multiple admissible control strategies coexist within the optimization landscape.

\begin{figure}[H]
\centering

\includegraphics[width=0.85\textwidth]{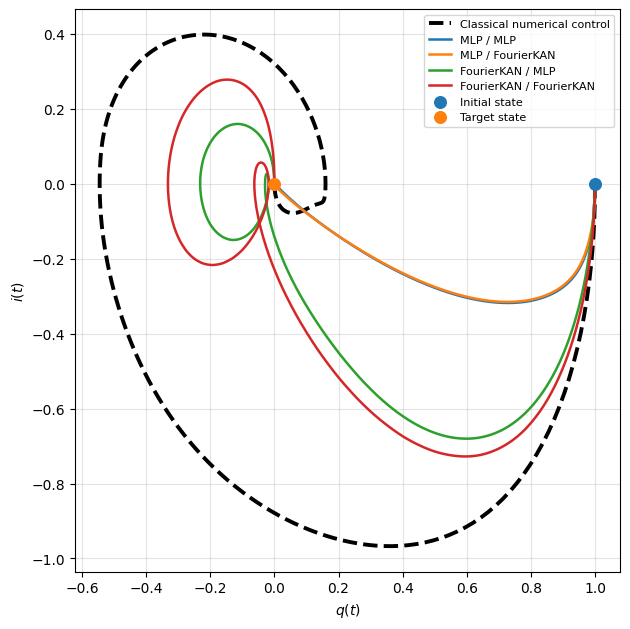}

\caption{
Phase portrait
}

\end{figure}

\subsection*{Spectral Analysis and Quantitative Metrics}

Figure 6 shows the Fourier spectra associated with the learned controls together with the classical nonlinear reference solution.

\begin{figure}[H]
\centering

\includegraphics[width=0.85\textwidth]{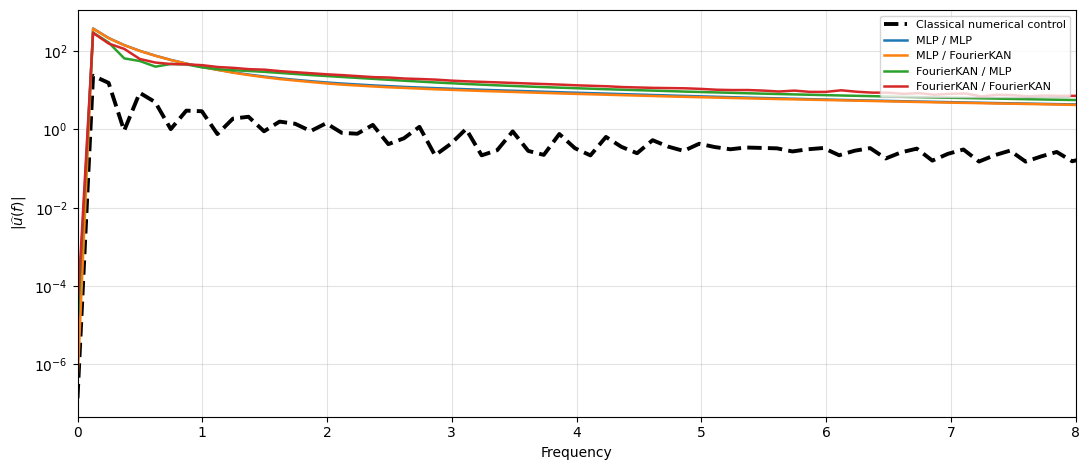}

\caption{
Fourier spectra 
}

\end{figure}

Compared with the classical optimized control, the neural architectures generate substantially different frequency distributions. As in the linear RLC problem, the Fourier-based configurations preserve stronger high-frequency components and exhibit richer oscillatory structures, whereas the MLP-based controls remain more concentrated in the low-frequency range.

These observations suggest that the frequency preferences associated with each architecture are not restricted to linear systems. Rather, they remain visible after the introduction of nonlinear dynamics, indicating that the architectural biases identified in the RLC experiments persist in a more complex setting.
Table 3 summarizes the main quantitative metrics associated with the nonlinear Duffing experiments. 

\begin{table}[H]
\centering
\small
\begin{tabular}{|l|c|c|c|c|c|}
\hline
Control Strategy & $E(u)$ & Smoothness & Spectral Centroid & $q(T)$ & $q'(T)$ \\
\hline

Classical nonlinear control
& 0.410647
& 1.005280
& --
& $-0.000200$
& $0.000100$ \\
\hline

MLP / MLP
& 1.625554
& 0.902770
& 13.422572
& $0.000100$
& $-0.000200$ \\
\hline

MLP / FourierKAN
& 1.630697
& 0.960536
& 13.610581
& $0.000200$
& $0.000100$ \\
\hline

FourierKAN / MLP
& 1.012200
& 1.305056
& 16.915188
& $-0.000300$
& $0.000200$ \\
\hline

FourierKAN / FourierKAN
& 1.005260
& 2.163774
& 17.735747
& $-0.000200$
& $0.000300$ \\
\hline

\end{tabular}
\caption{Summary of the main quantitative metrics associated with the different nonlinear control strategies.}
\label{tab:duffing_metrics}
\end{table}

The spectral-centroid values provide quantitative support for the trends observed in Figure 6. Architectures incorporating Fourier-based representations consistently exhibit larger spectral centroids, indicating a systematic shift toward higher frequencies. The same ordering observed in the linear RLC experiments is preserved in the nonlinear regime, suggesting that the spectral characteristics of the learned controls are strongly influenced by the underlying neural architecture.

At the same time, all configurations achieve very small terminal errors in both \(q(T)\) and \(q'(T)\). Consequently, the observed differences cannot be attributed to failures in solving the control problem. Instead, they reflect different admissible control strategies capable of reaching the same target state.

\subsection{Persistence of Architecture-Dependent Functional Priors}

The numerical experiments indicate that the architecture-dependent
effects observed in the linear RLC problem persist after the
introduction of nonlinear dynamics and, in some cases, appear to
become more pronounced.

Although all configurations achieve comparable terminal accuracies, they do so through substantially different control strategies. This behavior suggests that the underlying optimization problem admits multiple physically admissible solutions capable of satisfying the same dynamical objectives.

From this perspective, the neural architecture acts as an implicit functional prior that influences which region of the admissible-solution space is preferentially explored during training. Rather than serving solely as an approximation mechanism, the architecture biases the optimization process toward particular families of controls.

The persistence of these effects across both the linear RLC system and the nonlinear Duffing oscillator suggests that architecture-dependent functional specialization may represent a general phenomenon in physics-informed control problems rather than a property of a specific dynamical system.

%% file: sections/section5_conclusions.tex
The primary objective of this work was to investigate how neural architectures influence the qualitative nature of the controls learned within a Physics-Informed Neural Network framework \cite{raissi2019}, but also the functional characteristics of the controls learned in dynamical systems

To address this question, we studied two controlled systems of increasing complexity: a linear RLC circuit and a nonlinear Duffing oscillator. In both cases, the numerical experiments revealed that different neural architectures systematically generate different control strategies despite being trained under identical governing equations, optimization objectives, physical constraints, and training settings.

Several consistent patterns emerged across the experiments. First, architectures based on Fourier representations tend to generate controls with richer oscillatory structures and larger spectral centroids. In contrast, MLP-based architectures generally produce smoother controls with a stronger concentration of low-frequency components. These differences are reflected not only in the control signals themselves, but also in the geometry of the associated phase-space trajectories and in the energetic characteristics of the resulting solutions.

More importantly, the architecture-dependent effects observed in the linear RLC system persist after introducing nonlinear dynamics through the Duffing equation. The preservation of the same qualitative trends across both systems suggests that the phenomenon is not restricted to a particular model or dynamical regime. Instead, it appears to reflect a more general interaction between neural architectures and the optimization process underlying physics-informed control.

The observed behavior is also consistent with theoretical and empirical studies on spectral bias in deep neural networks \cite{rahaman2019spectral,xu2020frequency}. In particular, the MLP-based architectures considered here tend to generate smoother controls with a stronger concentration of low-frequency components, whereas Fourier-based architectures preserve richer oscillatory structures and larger spectral centroids. Although the present experiments do not establish a direct causal connection between spectral bias and control generation, the persistence of the same spectral ordering across both the linear RLC system and the nonlinear Duffing oscillator suggests that frequency-dependent learning effects play a significant role in shaping the learned controls.

Another notable observation is the emergence of a functional specialization phenomenon. Across both dynamical systems considered in this work, Fourier-based architectures appear naturally suited to representing oscillatory dynamical behavior, while MLP-based architectures tend to generate smoother and energetically more regular control signals. This observation is particularly interesting in light of recent studies comparing PINNs and KAN-based architectures \cite{liu2024kan,benchmark2026kanpinn, Somvanshi_2025}, which suggest that architectural choices may influence not only approximation accuracy but also the qualitative properties of the learned solutions.

Taken together, these observations support the interpretation of neural architectures as implicit functional priors. Rather than acting as neutral approximation tools, they influence which admissible solutions are preferentially selected during training. The resulting controls may therefore differ substantially even when they satisfy the same physical constraints and optimization objectives.

From the perspective of control theory \cite{zhou1996robust,brunton2022data}, these results suggest that architecture selection should be viewed not only as a computational choice but also as a modeling decision. If similar phenomena persist in more complex systems, the choice of neural architecture may influence not only how a control is computed, but also which physically admissible control is ultimately selected by the optimization process.

Several limitations should be acknowledged. The present study is restricted to relatively simple ordinary differential equation models and a limited set of neural architectures. In addition, the reported results correspond to representative numerical experiments rather than large-scale statistical analyses over multiple random initializations, optimization settings, and hyperparameter configurations. Consequently, the present work should be interpreted as evidence for the existence of architecture-dependent control-selection effects rather than as a complete characterization of their magnitude.

\subsection*{Future Work}

The results reported here open several directions for future research. A natural extension consists of investigating whether similar phenomena appear in higher-dimensional systems, PDE-constrained control problems, and quantum-control settings, where oscillatory behavior and spectral structure play a central role.

Another promising direction concerns the incorporation of spectral information directly into the learning process. Rather than relying exclusively on traditional energy and smoothness regularization terms, future control formulations could include explicit spectral constraints or frequency-domain objectives. Such approaches may provide finer control over the functional characteristics of the learned solutions.

From a methodological perspective, time-frequency representations such as spectrograms may offer a richer description of architecture-dependent effects than Fourier spectra alone. These representations could help characterize how different architectures distribute control effort across both time and frequency scales.

More broadly, a deeper theoretical understanding of neural architectures as implicit variational mechanisms remains an open problem. Establishing rigorous connections between spectral bias, implicit regularization, optimization dynamics, and control selection may provide a useful theoretical foundation for the design of future Scientific Machine Learning control frameworks.